\documentclass[12pt]{article}

\title{\textbf{The Halting Probability Omega: Irreducible Complexity in\\ Pure Mathematics\footnote
{Enriques lecture given Monday, October 30, 2006, at the University of Milan.}
}}

\author{\textbf{Gregory Chaitin}\thanks{IBM T. J. Watson Research Center, P. O. Box 218,
Yorktown Heights, NY 10598, U.S.A., \emph{chaitin@us.ibm.com.}}}

\date{}

\begin{document}

\maketitle

\begin{abstract}
Some G\"odel centenary reflections on whether incompleteness is really serious,
and whether mathematics should be done somewhat differently, based on using algorithmic
complexity measured in bits of information.
\end{abstract}
     
\section*{Introduction: What is mathematics?}
      
It is a pleasure for me to be here today giving this talk in a lecture series in honor
of Frederigo Enriques. Enriques was a great believer in mathematical intuition,
and disdained formal proofs. The work of G\"odel, Turing and myself that I will
review goes some way to justifying Enriques's belief in intuition.
And, as you will see, I also agree with Enriques's emphasis on the
importance of the philosophy and the history of science and mathematics.
    
This year is the centenary of Kurt G\"odel's birth. Nevertheless, his famous 1931
incompleteness theorem remains controversial.  To postmodernists, it justifies
the belief that truth is a social construct, not absolute.  Most mathematicians
ignore incompleteness, and carry on as before, in a formalist, axiomatic, Hilbertian, Bourbaki
spirit. I, on the contrary, have bet my life on the hunch that incompleteness is really serious,
that it cannot be ignored, and that it means that mathematics is actually somewhat different from
what most people think it is.
    
G\"odel himself did not think that his theorem showed that mathematics has limitations.
In several essays he made it clear that he believed that mathematicians could eventually
settle any significant question by using their mathematical intuition, their ability to
directly perceive the Platonic world of mathematical ideas, and by inventing or discovering new concepts
and new axioms, new principles. 
   
Furthermore, I share Enriques's faith in intuition. I think that excessive formalism and
abstraction is killing mathematics. In my opinion math papers shouldn't attempt to replace all words by formulas,
instead they should be like literary essays, they should attempt to explain and convince.
    
So let me tell you the story of metamathematics, of how mathematicians have tried to use mathematical
methods to study the power and the limitations of math itself. It's a fairly dramatic story; in a previous
era it might have been the subject of epic poems, of Iliads and Odysseys of verse.
I'll start with David Hilbert.
     
\section*{Hilbert: Can mathematics be entombed in a formal axiomatic theory?}
    
Hilbert stated the traditional belief that mathematics can provide absolute truth, complete certainty,
that mathematical truth is black or white with no uncertainty.
His contribution was to realize, to emphasize, that if this were the case, then there should be, there
ought to be, a formal axiomatic theory, a theory of everything, for all of mathematics.
    
In practice, the closest we have come to this today is Zermelo-Fraenkel set theory with the axiom of choice,
the formal theory ZFC using first-order logic, which seems to suffice for most contemporary mathematics.
      
Hilbert did not invent mathematical logic, he took advantage of work going back to Leibniz, de Morgan, Boole,
Frege, Peano, Russell and Whitehead, etc.
But in my opinion he enunciated more clearly than anyone before him the idea that if math provides absolute truth,
complete certainty, then there should be a finite set of axioms that we can all agree on from which it would
in principle be possible to prove all mathematical truths by mechanically following the rules of formal mathematical
logic.  It would be slow, but it would work like an army of reason marching inexorably forward.
It would make math into a merciless machine.
    
Hilbert did not say that mathematics should actually be done in this extremely formal way
in which proofs are broken down into their atomic steps, with nothing omitted, excruciatingly detailed, 
using symbolic logic instead of a normal human language.  But the idea was to eliminate 
all uncertainty, to make clear exactly when a proof is valid, so that this can be checked mechanically,
thus making mathematical truth completely objective, eliminating all subjective elements, all matters of opinion.
    
Hilbert started with Peano arithmetic, but his ultimate goal was to axiomatize analysis and then all of mathematics,
absolutely everything.
In 1931, however, G\"odel surprised everyone by showing that it couldn't be done, it was impossible.
   
\section*{G\"odel: ``This statement is unprovable!''}
      
In fact, G\"odel showed that no finite set of axioms suffice for elementary number theory, for the theory
of 0, 1, 2, \ldots\ and addition and multiplication, that is, for Peano arithmetic.
His proof is very strange. First of all he numbers all possible assertions and all possible proofs in
Peano arithmetic. This converts the assertion that $x$ is a proof of $y$ 
into an arithmetic assertion about $x$ and $y$.
    
Next G\"odel constructs an assertion that refers to itself indirectly. It says that if you calculate
a certain number, that gives you the number of an unprovable assertion, and this is done in such a way that
we get an arithmetic statement asserting its own unprovability.
    
Consider ``I am unprovable.'' It is either provable or not. If provable, we are proving a false assertion,
which we very much hope is impossible.  The only alternative left is that ``I'm unprovable'' is unprovable.
If so it is true but unprovable, and there is a hole in formal mathematics, a true assertion that we cannot prove.
In other words, our formal axiomatic theory must be incomplete, if we assume that only true assertions can be proved,
which we fervently hope to be the case. Proving false assertions is even worse than not being able to prove
a true assertion!
    
So that's G\"odel's famous 1931 incompleteness theorem, and it was a tremendous shock to everyone.
When I was a young student I read essays by John von Neumann, Hermann Weyl and others attesting to what a shock it was.
I didn't realize that the generation that cared about this had been swept away by the Second World War and that
mathematicians were going on exactly as before, ignoring G\"odel.  In fact, I thought that what G\"odel
discovered was only the tip of the iceberg. I thought that the problem had to be really serious, really profound, and
that the traditional philosophy of math couldn't be slightly wrong, so that even a small scratch would shatter it into 
pieces. It had to be all wrong, in my opinion.
    
Does that mean that if you cannot prove a result that you like and have some numerical evidence for, if you cannot
do this in a week, then invoking G\"odel, you just add this conjectured result as a new axiom?!
No, not at all, that is too extreme a reaction. But, as I will explain later, I do have something like that in mind.
    
You see, the real problem with G\"odel's proof is that it gives no idea how serious incompleteness is.
G\"odel's true but unprovable assertion is bizarre, so it is easy to shrug it off.
But if it turns out that incompleteness is pervasive, is ubiquitous, that is another matter.
    
And an important first step in the direction of showing that incompleteness is really serious was taken
only five years later, in 1936, by Alan Turing, in a famous paper ``On computable numbers\ldots''
    
\section*{Turing: Most real numbers are uncomputable!}
      
This paper is remembered, in fact, celebrated nowadays, for proposing a toy model of the computer called a Turing machine,
and for its discussion of what we now call the halting problem.
But this paper is actually about distinguishing between computable and uncomputable real numbers,
numbers like $\pi$ or $e$ or $\sqrt{2}$ that we can compute with infinite precision, with arbitrary accuracy,
and those that we cannot.
      
Yes, it's true, Turing's paper does contain the idea of software as opposed to hardware, of a universal digital
machine that can simulate any other special-purpose digital machine.  Mathematicians refer to this as a universal Turing
machine, and, as I learned from von Neumann, it is the conceptual basis for all computer technology.
But even more interesting is the fact that it is easy to see that most real numbers are uncomputable, and
the new perspective this gives on incompleteness, as Turing himself points out.
Let me summarize his discussion.
    
First of all, all possible software, all possible algorithms, can be placed in an infinite list and numbered
1, 2, 3, \ldots\ and so this set is denumerable or countable.
However, the set of real numbers is, as Cantor showed, a higher-order infinity, it's uncountable, nondenumerable.
Therefore most real numbers must be uncomputable.  In Turing's paper he exhibits a single example of an uncomputable
real, one that is obtained by applying Cantor's diagonal method to the list of all computable real numbers to obtain
a new and different real, one that is not in the list.
    
In fact, using ideas that go back to Emile Borel, it is easy to see that if you pick a real number between 0 and 1
at random, with uniform probability, it is possible to pick a computable real, but this is infinitely improbable,
because the computable reals are a set of measure zero, they can be covered with intervals whose total length is
arbitrarily small. Just cover the first computable real in [0,1) with an interval of size $\epsilon/2$, the second
real with an interval of size $\epsilon/4$, and in general the $N$th real with an interval of size
$\epsilon/2^N$. This covering has lengths totalling exactly $\epsilon$, which can be made as small as we want.
    
Turing does not, however, make this observation. Instead he points out that it looks easy to compute his uncomputable
real by taking the $N$th digit produced by the $N$th program and changing it. Why doesn't this work? Because we can never
decide if the $N$th program will ever produce an $N$th digit! If we could, we could actually diagonalize over all
computable reals and calculate an uncomputable real, which is impossible.
And being able to decide if the $N$th program will ever output an $N$th digit is a special case of Turing's famous
halting problem.
    
Note that there is no problem if you place an upper bound on the time allowed for a computation.
It is easy to decide if the $N$th program outputs an $N$th digit in a trillion years, all you need to do is be patient and
try it and see.  Turing's halting problem is only a problem if there is no time limit.
In other words, this is a deep conceptual problem, not a practical limitation on what we can do.
    
And, as Turing himself points out, incompleteness is an immediate corollary. For let's say we'd like to be able
to prove whether individual computer programs, those that are self-contained and read no input, eventually halt or not.
There can be no formal axiomatic theory for this, because if there were, by systematically running through
the tree of all possible proofs, all possible deductions from the axioms using formal logic, we could always
eventually decide whether an individual program halts or not, which is impossible.
     
In my opinion this is a fundamental step forward in the philosophy of mathematics because it makes incompleteness
seem much more concrete and much more natural. It's almost a problem in physics, it's about a machine, you just
ask whether or not it's going to eventually stop, and it turns out there's no way, no general way, to answer that 
question.
     
Let me emphasize that if a program does halt, we can eventually discover that. The problem, an extremely deep one,
is to show that a program will never halt if this is in fact so.  One can settle many special cases, even
an infinity of them, but no finite set of axioms can enable you to settle all possible cases.
    
My own work takes off from here. My approach to incompleteness follows Turing, not G\"odel. Later I'll consider
the halting probability $\Omega$ and show that this number is wildly, in fact maximally, uncomputable and unknowable.
I'll take Turing's halting problem and convert it into a real number\ldots
   
My approach is very 1930's. All I add to Turing is that I measure software complexity, I look at the size of computer programs.
In a moment, I'll tell you how these ideas actually go back to Leibniz.
But first, let me tell you more about Borel's ideas on uncomputable reals, which are closely related to Turing's ideas.
    
\section*{Borel: Know-it-all and unnameable reals}
      
Borel in a sense anticipated Turing, because he came up with an example of an uncomputable real in a small paper published in 1927.
Borel's idea was to use the digits of a single real number as an oracle that can answer any yes/no question.
Just imagine a list of all possible yes/no questions in French, said Borel. This is obviously a countable infinity,
so there is an $N$th question, an $N+1$th question, etc.   And you can place all the answers in the decimal expansion of
a single real number; just use the $N$th digit to answer the $N$th question.
Questions about history, about math, about the stock market!
    
So, says Borel, this real number exists, but to him it is a mathematical fantasy, not something real.
Basically Borel has a constructive attitude, he believes that something exists only if we can calculate it,
and Borel's oracle number can certainly not be calculated.
    
Borel didn't linger over this, he made his point and moved on, but his example is in my opinion a very interesting one,
and will later lead us step by step to my $\Omega$ number.
     
Before I used Borel's ideas on measure and probability to point out that Turing's computable reals have measure zero,
they're infinitely unlikely.  Borel however seemed unaware of Turing's work. His own version of these ideas, in his
final book, written when he was 80, \emph{Les nombres inaccessibles,} is to point out that the set of reals that can
somehow be individually identified, constructively or not, has measure zero, because the set of all possible
descriptions of a real is countable.  Thus, with probability one, a real cannot be uniquely specified, it can never
be named, there are simply not enough names to go around!
    
The real numbers are the simplest thing in the world geometrically, they are just points on a line.
But arithmetically, as individuals, real numbers are actually rather unreal.  Turing's 1936 uncomputable real is
just the tip of the iceberg, the problem is a lot more serious than that.
     
Let me now talk about looking at the size of computer programs and what that has to tell us about incompleteness.
To explain why program size is important, I have to start with Leibniz, with some ideas in his
1686 \emph{Discours de m\'etaphysique,} which was found among his papers long after his death.
    
\section*{Theories as software, Understanding as compression, Lawless incompressible facts}
      
The basic model of what I call algorithmic information theory (AIT) is that a scientific theory is a computer
program that enables you to compute or explain your experimental data:
\begin{center}
theory (program) $\longrightarrow$ \textbf{Computer} $\longrightarrow$ data (output).
\end{center}
In other words, the purpose of a theory is to compute facts.
The key observation of Leibniz is that there is always a theory that is as complicated as the facts it is
trying to explain. This is useless: a theory is of value only to the extent that it compresses a great many
bits of data into a much smaller number of bits of theory.
    
In other words, as Hermann Weyl put it in 1932, the concept of law becomes vacuous if an arbitrarily complicated
law is permitted, for then there is always a law.  A law of nature has to be much simpler than the data it explains,
otherwise it explains nothing.  The problem, asks Weyl, is how can we measure complexity?
Looking at the size of equations is not very satisfactory.
     
AIT does this by considering both theories and data to be digital information; both are a finite string of bits.
Then it is easy to compare the size of the theory with the size of the data it supposedly explains, by merely
counting the number of bits of information in the software for the theory and comparing this with the number of
bits of experimental data that we are trying to understand.
    
Leibniz was actually trying to distinguish between a lawless world and one that is governed by law. He was trying
to elucidate what it means to say that science works. This was at a time when modern science, then called
mechanical philosophy, was just beginning; 1686 was the year before Leibniz's nemesis Newton published his \emph{Principia.}
    
Leibniz's original formulation of these ideas was like this. Take a piece of paper, and spot it with a quill pen, so
that you get a finite number of random points on a page. There is always a mathematical equation that passes precisely
through these points. So this cannot enable you to distinguish between points that are chosen at random and points
that obey a law.  But if the equation is simple, then that's a law.
If, on the contrary, there is no simple equation, then the points are lawless, random.
     
So part and parcel of these ideas is a definition of randomness or lawlessness for finite binary strings, as those
which cannot be compressed into a program for calculating them that is substantially smaller in size.
In fact, it is easy to see that most finite binary strings require programs of about the same size as they are.
So these are the lawless, random or algorithmically irreducible strings, and they are the vast majority of all strings.
Obeying a law is the exception, just as being able to name an individual real is an exception.
    
Let's go a bit further with this theories as software model.  Clearly, the best theory is the simplest, the most concise, 
the smallest
program that calculates your data. So let's abstract things a bit, and consider what I call \textbf{elegant} programs:
\begin{itemize}
\item
A program is elegant if no smaller program written in the same language, produces the same output.
\end{itemize}
In other words, an elegant program is the optimal, the simplest theory for its output.
How can we be sure that we have the best theory?
How can we tell whether a program is elegant?
The answer, surprisingly enough, is that we can't!
    
\section*{Provably elegant programs}
       
To show this, consider the following paradoxical program $P$:
\begin{itemize}
\item
$P$ computes the output of the first provably elegant program larger than $P$.
\end{itemize}
In other words, $P$ systematically deduces all the consequences of the axioms,
which are all the theorems in our formal axiomatic theory.  As it proves each theorem, $P$ examines
it.  First of all, $P$ filters out all proofs that do not show that a particular program
is elegant. For example, if $P$ finds a proof of the Riemann hypothesis, it throws that away;
it only keeps proofs that programs are elegant.  And as it proves that individual programs are elegant, 
it checks each provably elegant program to see if this program is larger than $P$.
As soon as $P$ finds a provably elegant program that is larger than it is, it
starts running that program, and produces that program's output as its own output. In other words,
$P$'s output is precisely the same as the output of the first provably elegant program
that is larger than $P$.
    
However, $P$ is too small to produce the same output as an elegant program that is 
larger than $P$, because this contradicts the definition of elegance!
What to do? How can we avoid the contradiction?
    
First of all, we are assuming that our formal axiomatic theory only proves true theorems,
and in particular, that if it proves that a program is elegant, this is in fact the case.
Furthermore, the program $P$ is not difficult to write out; I've done this in one
of my books using the programming language called LISP.  So the only way out is if $P$
never finds the program it is looking for!
In other words, the only way out is if it is never possible to prove that a program that's
larger than $P$ is elegant!
But there are infinitely many possible elegant programs, and they can be arbitrarily big.
But provably elegant programs can't be arbitrarily big, they can't be larger than $P$.
     
So how large is $P$, that's the key question. Well, the bulk of $P$ is actually
concerned with systematically producing all the theorems in our formal axiomatic theory.
So I'll define the complexity of a formal axiomatic theory to be the size in bits of the
smallest program for doing that.
Then we can restate our metatheorem like this: You can't prove that a program is elegant
if its size in bits is substantially larger than the complexity of the formal axiomatic
theory you are working with.
In other words, using a formal axiomatic theory with $N$ bits of complexity,
you can't prove that any program larger than $N + c$ bits in size is elegant.
Here the constant $c$ is the size in bits of the main program in $P$, the fixed
number of bits not in that big $N$-bit subroutine for running the formal axiomatic
theory and producing all its theorems.
    
Loosely put: 
\begin{itemize}
\item
You need an $N$-bit theory to show that an $N$-bit program is elegant.\footnote
{By the way, it is an immediate corollary that the halting problem is unsolvable,
because if we could decide which programs halt, then we could run all the programs that halt
and see what they output, 
and this would give us a way to determine which halting programs are elegant, which we've just shown is impossible.
This is a new information-theoretic proof of Turing's theorem, rather different from Turing's original
diagonal-argument proof.}
\end{itemize}
   
Why is this so interesting?
Well, right away it presents incompleteness in an entirely new light.
How?
Because it shows that mathematics has infinite complexity, but any formal axiomatic theory
can only capture a finite part of this complexity. In fact, just knowing which programs are elegant
has infinite complexity. 
    
So this makes incompleteness very natural; math has an infinite basis,
no finite basis will do.
Now incompleteness is the most natural thing in the world; it's not at all mysterious!
    
Now let me tell you about the halting probability $\Omega$, which shows even better that math is
infinitely complex.
    
\section*{What is the halting probability $\Omega$?}
     
Let's start with Borel's know-it-all number, but now let's use the $N$th binary digit to tell us whether or
not the $N$th computer program ever halts. So now Borel's number is an oracle for the halting problem.
For example, there is a bit which tells us whether or not the Riemann hypothesis is true, for that
is equivalent to the statement that a program that systematically searches for zeros of the zeta function
that are in the wrong place, never halts.
     
It turns out that this number, which I'll call Turing's number even though it does not occur in Turing's paper,
is wasting bits, it is actually highly redundant.  We don't really need $N$ bits to answer $N$ cases of the
halting problem, a much smaller number of bits will do. 
Why?
     
Well, consider some large number $N$ of cases of the halting problem, some large number $N$ of individual
programs for which we want to know whether or not each one halts.  Is this really $N$ bits of mathematical information?
No, the answers are not independent,  they are highly correlated. How?
Well, in order to answer $N$ cases of the halting problem, we don't really need to know each individual answer;
it suffices to know how many of these $N$ programs will eventually halt. Once we know this number, which is
only about $\log_2 N$ bits of information, we can run the $N$ programs in parallel until exactly this number
of them halt, and then we know that none of the remaining programs will ever halt.
And $\log_2 N$ is much, much less than $N$ for all sufficiently large $N$.
In other words, Turing's number isn't the best possible oracle for the halting problem. It is highly redundant,
it uses far too many bits.
     
Using essentially this idea, we can get the best possible oracle number for the halting problem;
that is the halting probability $\Omega$, which has no redundancy, none at all.
    
I don't have time to explain this in detail, but here is a formula for the halting probability:
\[
\Omega = \sum_{\mbox{\scriptsize $p$ halts}} 2^{-|p|}.
\]
The idea is that each $K$-bit program that halts contributes exactly $1/2^K$ to
the halting probability $\Omega$.
In other words, $\Omega$ is the halting probability of a program $p$ whose bits are generated by independent
tosses of a fair coin. 
    
\textbf{Technical point:} For this to work, for this sum to converge to a number between 0 and 1 instead of
diverging to infinity, it is important that programs be self-delimiting, that no extension of a valid
program be a valid program.  In other words, our computer must decide by itself when to stop reading
the bits of the program without waiting for a blank endmarker. In old Shannon-style information theory
this is a well-known lemma called the Kraft inequality that applies to prefix-free sets of strings,
to sets of strings which are never prefixes or extensions of each other.
And an extended, slightly more complicated, version of the Kraft inequality plays a fundamental role in AIT.
    
I should also point out that the precise numerical value of $\Omega$ depends on your choice of computer
programming language, or, equivalently, on your choice of universal self-delimiting Turing machine.
But its surprising properties do not, they hold for a large class of universal Turing machines.
     
Anyway, once you fix the programming language, the precise numerical value of $\Omega$ is determined, it's well-defined.
Let's imagine having $\Omega$ written out in base-two binary notation:
\[
\Omega = .110110\ldots
\]
These bits are totally lawless, algorithmically irreducible mathematical facts.  They cannot be compressed into
any theory smaller than they are. 
    
More precisely, the bits of the halting probability $\Omega$ are both computationally and logically irreducible:
\begin{itemize}
\item
You need an $N$-bit program to calculate $N$ bits of $\Omega$ (any $N$ bits, not just the first $N$).
\item
You need an $N$-bit theory to be able to determine $N$ bits of $\Omega$ (any $N$ bits, not just the first $N$).
\end{itemize}
      
$\Omega$ is an extreme case of total lawlessness; in effect, it shows that God plays dice in pure mathematics.
More precisely, the bits of $\Omega$ refute Leibniz's principle of sufficient reason, because they are
mathematical facts that are true for no reason (no reason simpler than they are). Essentially the only way
to determine bits of $\Omega$ is to directly add these bits to your axioms. But you can prove anything by
adding it as a new axiom; that's not using reasoning!
    
Why does $\Omega$ have these remarkable properties?
Well, because it's such a good oracle for the halting problem.
In fact, knowing the first $N$ bits of $\Omega$ enables you to answer the halting problem for all
programs up to $N$ bits in size. 
And you can't do any better; that's why these bits are incompressible, irreducible information.
If you think of what I called Turing's number as a piece of coal, then $\Omega$ is the diamond that you get from this
coal by subjecting it to very high temperatures and pressures.
A relatively small number of bits of $\Omega$ would in principle enable you to tell whether or not the Riemann hypothesis is false.
     
\section*{Concluding discussion}
      
So $\Omega$ shows us, directly and immediately, that math has infinite complexity, because the bits of $\Omega$ are
infinitely complex.
But any formal axiomatic theory only has a finite, in fact, a rather small complexity, otherwise we wouldn't believe in it!
What to do? How can we get around this obstacle?
Well, by increasing the complexity of our theories, by adding new axioms, complicated axioms that are pragmatically
justified by their usefulness instead of simple self-evident axioms of the traditional kind.
     
Here are some recent examples:
\begin{itemize}
\item
The hypothesis that
\textbf{P} $\neq$ \textbf{NP} in theoretical computer science.
\item
The axiom of projective determinacy in abstract set theory.
\item
Various versions of the Riemann hypothesis in analytic number theory.
\end{itemize}
     
In other words, I am advocating a ``quasi-empirical'' view of mathematics, a term that was invented by Imre Lakatos, by the way.
(He wouldn't necessarily approve of the way I'm using it, though.)
     
To put it bluntly, from the point of view of AIT, mathematics and physics are not that different.
In both cases, theories are compressions of facts, in one case facts we discover in a physics lab, in the other case,
numerical facts discovered using a computer. Or, as Vladimir Arnold so nicely puts it, math is like physics, except that
the experiments are cheaper!
I'm not saying that math and physics are the same, but I am saying that maybe they are not as different as most people think.
     
Another way to put all of this, is that the DNA for pure math, $\Omega$, is infinitely complex, whereas the human genome
is $3 \times 10^9$ bases = $6 \times 10^9$ bits, a large number, but a finite one. So pure math is even
more complex than the traditional domain of the complicated, biology!
Math does not have finite complexity the way that Hilbert thought, not at all, on the contrary!
    
These are highly heretical suggestions, suggestions that the mathematics community is extremely uncomfortable with.
And I have to confess that I have attempted to show that incompleteness is serious and that math should be done somewhat
differently, but I haven't been able to make an absolutely watertight case. I've done the best I can with one lifetime 
of effort, though.
     
But if you really started complicating mathematics by adding new non-self-evident axioms, what would happen?
Might mathematics break into separate factions? Might different groups with contradictory axioms go to war?
Hilbert thought math was an army of reason marching inexorably forward, but this sounds more like anarchy!
    
Perhaps anarchy isn't so bad; it's better than a prison, and it leaves more room for intuition and creativity.
I think that Enriques might have been sympathetic to this point of view. After all, as Cantor, who created
a crazy, theological, paradoxical theory of infinite magnitudes, said, the essence of mathematics resides in its freedom,
in the freedom to imagine and to create.
    
\section*{Bibliography}

For more on this, please see my book [4] that has just been published in Italian, or my previous book,
which is currently available in three separate editions [1, 2].
Some related books in Italian [3, 5]  and papers [6--8] are also listed below.
\begin{itemize}
\item[{[1]}]
G. Chaitin, \emph{Meta Math!,} Pantheon, New York, 2005 (hardcover), Vintage, New York, 2006 (softcover).
\item[{[2]}]
G. Chaitin, \emph{Meta Maths,} Atlantic Books, London, 2006.
\item[{[3]}]
U. Pagallo, \emph{Introduzione alla filosofia digitale. Da Leibniz a Chaitin,} Giappichelli, Turin, 2005.
\item[{[4]}]
G. Chaitin, \emph{Teoria algoritmica della complessit\`a,} Giappichelli, Turin, 2006.
\item[{[5]}]
U. Pagallo, \emph{Teoria giuridica della complessit\`a,} Giappichelli, Turin, 2006.
\item[{[6]}]
G. Chaitin, ``How real are real numbers?,'' \emph{International Journal of Bifurcation and Chaos} 
\textbf{16} (2006), pp.\ 1841--1848.
\item[{[7]}]
G. Chaitin, ``Epistemology as information theory: from Leibniz to $\Omega$,'' \emph{Collapse} \textbf{1} (2006), pp.\ 27--51.
\item[{[8]}]
G. Chaitin, ``The limits of reason,'' \emph{Scientific American} \textbf{294}, No.\ 3 (March 2006), pp.\ 74--81.
(Also published as ``I limiti della ragione,'' \emph{Le Scienze,} May 2006, pp.\ 66--73.)
\end{itemize}

\end{document}